\documentclass{amsart}
\usepackage{amsmath, amssymb}
\usepackage{amsfonts}
\usepackage[arrow,matrix,curve,cmtip,ps]{xy}

\usepackage{amsthm}

\allowdisplaybreaks

\newtheorem{theorem}{Theorem}[section]
\newtheorem{lemma}[theorem]{Lemma}
\newtheorem{proposition}[theorem]{Proposition}
\newtheorem{corollary}[theorem]{Corollary}

\newtheorem*{theorem*}{Theorem}
\theoremstyle{remark}
\newtheorem{remark}[theorem]{Remark}
\newtheorem{definition}[theorem]{Definition}
\newtheorem{example}[theorem]{Example}

\numberwithin{equation}{section}

\newcommand{\T}{\mathbb{T}}

\newcommand{\B}{\mathcal{B}}

\newcommand{\M}{\mathcal{M}}
\newcommand{\OA}{\mathcal{O}_A}
\newcommand{\Gr}{\mathrm{Gr}}

\newcommand{\G}{\mathcal{G}}
\newcommand{\F}{\mathcal{F}}

\newcommand{\Pow}{\mathcal{P}}
\newcommand{\Hi}{\mathcal{H}}

\newcommand{\aut}{\operatorname{Aut}}

\begin{document}
\title[Exel-Laca algebras and $C^*$-algebras
associated to graphs]{A unified approach to Exel-Laca
algebras and \boldmath{$C^*$}-algebras associated to graphs}

\author{Mark Tomforde 
}

\address{Department of Mathematics\\ Dartmouth College\\
Hanover\\ NH 03755-3551\\ USA}
\email{mark.tomforde@dartmouth.edu}


\date{\today}

\begin{abstract}

We define an ultragraph, which is a generalization of a
directed graph, and describe how to associate a
$C^*$-algebra to it.  We show that the class of ultragraph
algebras contains the $C^*$-algebras of graphs as well as
the Exel-Laca algebras.  We also show that many of the
techniques used for graph algebras can be applied to ultragraph algebras and that the ultragraph provides a useful
tool for analyzing Exel-Laca algebras.  Our results include
versions of the Cuntz-Krieger Uniqueness Theorem and the
Gauge-Invariant Uniqueness Theorem for ultragraph
algebras.

\end{abstract}

\maketitle

\section{Introduction}

In the early 1980's Cuntz and Krieger considered a class of
$C^*$-algebras that arose in the study of topological Markov
chains \cite{CK}.  These Cuntz-Krieger algebras $\OA$ are
generated by partial isometries whose relations are
determined by a finite matrix $A$ with entries in $\{ 0 , 1
\}$.  In order for their $C^*$-algebras to be unique, Cuntz
and Krieger assumed that the matrix $A$ also satisfied a
nondegeneracy condition called Condition~(I).  Since their
introduction Cuntz-Krieger algebras have been generalized
in a myriad of ways.  Two important
generalizations are the Exel-Laca algebras of \cite{EL} and
the $C^*$-algebras of directed graphs \cite{KPRR, KPR, BPRS,
FLR}.

The Exel-Laca algebras are in some sense the most direct
generalization of Cuntz-Krieger algebras.  In 1999 Exel and
Laca extended the definition of $\OA$ to allow for
infinite matrices \cite{EL}.  Furthermore, the only
restriction placed on these matrices was that they had no
zero rows.  Motivated by the $C^*$-algebras associated to
graphs, Exel and Laca avoided the need for Condition~(I) by
instead requiring the generating partial isometries to be
universal for their defining relations.  Since
Condition~(I) was imposed by Cuntz and Krieger to insure
uniqueness, this universal definition agreed with Cuntz and
Krieger's for finite matrices satisfying Condition~(I). 

The generalization of Cuntz-Krieger algebras to $C^*$-algebras
of directed graphs is slightly less direct.  In 1982 Watatani
noted that one could view $\OA$ as the $C^*$-algebra of a
finite directed graph with vertex adjacency matrix $A$
\cite{Wat}.  The fact that $A$ satisfied Condition~(I)
implied, among other things, that this graph had no sinks or
sources.  It was to be approximately 15 years, however,
before these graph ideas were explored more fully.  In the
late 1990's generalizations of these $C^*$-algebras were
considered for possibly infinite graphs that were allowed to
contain sinks and sources.  Originally, a definition was
given only for graphs that are row-finite; that is, each
vertex is the source of only finitely many edges \cite{KPRR,
KPR, BPRS}.  However, in 2000 this definition was extended
to obtain an appropriate notion of $C^*$-algebras for
arbitrary graphs
\cite{FLR}.

The relationship between Exel-Laca algebras and graph algebras
is somewhat subtle.  As mentioned before, the class
of Exel-Laca algebras of finite matrices satisfying
Condition~(I) and the class of $C^*$-algebras of
certain finite graphs both coincide with the
Cuntz-Krieger algebras.  However, the infinite case is more
complicated.  

Let $G$ be a graph with no sinks or sources.  The edge matrix
$A$ for the graph $G$ is the matrix indexed by the
edges of $G$ with $A(e,f)=1$ if $r(e) = s(f)$, and
$A(e,f)=0$ otherwise.  It was shown in \cite{FLR} that
$C^*(G)$ is canonically isomorphic to $\OA$.  Thus
$C^*$-algebras of graphs without sinks or sources are
Exel-Laca algebras.

Unfortunately, the reverse inclusion is not true.  It was
shown in \cite{RS} that there exist Exel-Laca algebras
that are not graph algebras.  However, one can obtain a
partial converse in the row-finite case.  If $A$ is a matrix,
then one may form a graph $\Gr (A)$ by defining the vertex
set of $\Gr (A)$ to be the index set of $A$, and defining
the number of edges from $v$ to $w$ to be $A(v,w)$.  If $A$
is a $\{0,1\}$-matrix that is row-finite (i.e.,
each row of $A$ is eventually zero), then $\Gr (A)$ is a
row-finite graph and $\OA \cong C^* ( \Gr (A))$.  This
isomorphism is obtained through the use of the dual graph of
$\Gr (A)$ and it is not canonical.  It was shown in
\cite{RS} that $C^*(\Gr(A))$ is always a
$C^*$-subalgebra of $\OA$, but when $A$ is not row-finite
this subalgebra may be very different from $\OA$ (see
\cite[Remark 15]{Szy} and \cite{DT2}).

These relationships are summarized in the following diagram:

\text{ }

$$\begin{matrix} \text{CK} & \subset & \text{F}' &
\subset & \text{RF}' & \subset &
\text{G}' & \subset & \text{EL} \\ 
& & \cap &  & \cap &  & \cap & &  \\
& & \text{F} & \subset & \text{RF} & \subset & \text{G} &  & 
\end{matrix}
$$

\text{}

\begin{align*} \text{CK} & = \text{Cuntz-Krieger algebras
$\OA$ with $A$ satisfying Condition~(I)}
\\ \text{F}' & = \text{$C^*$-algebras of finite graphs with
no sinks or sources} \\
\text{RF}' & = \text{$C^*$-algebras of
row-finite graphs with no sinks or sources} \\ 
& = \text{Exel-Laca algebras of row-finite matrices with no
zero rows or columns} \\
\text{G}' & = \text{$C^*$-algebras of
graphs with no sinks or sources} \\
\text{EL} & = \text{Exel-Laca algebras} \\
\text{F} & = \text{$C^*$-algebras of finite graphs} \\
\text{RF} & = \text{$C^*$-algebras of row-finite graphs} \\
\text{G} & = \text{$C^*$-algebras of graphs} 
\end{align*}

In some sense, it is unfortunate that there are Exel-Laca
algebras that are not graph algebras.  The graph $G$ is an
extremely useful tool for analyzing $C^*(G)$, and many
results take a more elegant form when stated in terms of
graphs rather than matrices.  In fact, when studying the
$\OA$'s Exel and Laca found it convenient to state
many of their hypotheses and results in terms of the graph
$\Gr (A)$.

In this paper we describe a generalized notion of graph,
which we call an ultragraph, and describe a way to
associate a $C^*$-algebra to it.  We shall see that
the $C^*$-algebras of ultragraphs with no
sinks and in which every vertex emits finitely many edges are
precisely the Exel-Laca algebras.  Thus the ultragraph
algebras fit into our diagram as follows:

\text{ }

$$\begin{matrix} \text{CK} & \subset & \text{F}' &
\subset & \text{RF}' & \subset &
\text{G}' & \subset & \text{EL}=\widetilde{\text{U}} \\ 
& & \cap &  & \cap &  & \cap & & \cap \\
& & \text{F} & \subset & \text{RF} & \subset & \text{G} &
\subset & \text{U}
\end{matrix}
$$

\begin{align*}
\widetilde{\text{U}} = & \ \text{$C^*$-algebras of
ultragraphs with no sinks and in} \\ & \quad \quad \quad 
\text{which every vertex emits finitely many edges.} \\ 
\text{U} = & \ \text{$C^*$-algebras of ultragraphs}
\\
\end{align*}

Therefore ultragraph algebras give an alternative
(and in the author's opinion, more convenient) way to view
Exel-Laca algebras.  In addition, they provide a
reasonable notion of ``Exel-Laca algebras with sinks". 
We shall see that for a countably indexed matrix $A$, one
may create an ultragraph $\G_A$ with edge matrix $A$
and for which $C^*(\G_A)$ is canonically isomorphic to
$\OA$.  Furthermore, the ultragraph $\G$
provides a useful tool for analyzing the structure of
$C^*(\G)$, just as the graph does for graph algebras.  

We shall prove in \S\ref{EL-sec} that the
$C^*$-algebras of ultragraphs with no sinks and no
infinite emitters are precisely the Exel-Laca algebras. 
However, when deducing results about ultragraph algebras
we will almost always prove them from first principals
rather than simply translating the corresponding result from
$\OA$ to the ultragraph setting.  This keeps our
treatment more self-contained, and more importantly it shows
that ultragraph algebras are tractable objects of study. 
In addition, since many of our techniques generalize those
used for graph algebras, this approach shows that by viewing
Exel-Laca algebras as ultragraph algebras one can forget
about matrix techniques entirely and instead apply the (often
more straightforward) graph techniques.

This paper is organized as follows.  In \S2 we define
an ultragraph $\G$ and describe a way to associate a
$C^*$-algebra $C^*(\G)$ to it.  In \S3 we discuss a
natural way to view graph algebras as ultragraph algebras
and discuss loops in ultragraphs.  In \S4 we show that the
$C^*$-algebras of ultragraphs with no singular vertices
are precisely the Exel-Laca algebras.  In \S5 we describe
a method for realizing certain subalgebras of $C^*(\G)$ as
graph algebras.  In \S6 we discuss methods to remove singular
vertices from ultragraphs by adding tails, and we obtain
versions of the Gauge-Invariant Uniqueness Theorem and the
Cuntz-Krieger Uniqueness Theorem for ultragraph algebras. 

The author would like to thank Doug Drinen for many valuable
discussions relating to this work.

\section{Ultragraphs and their $C^*$-algebras}

To provide motivation for the definition of an ultragraph,
recall that when $G$ is a graph with no sinks or sources
and with edge matrix $A$, then
$C^*(G)$ is canonically isomorphic to $\OA$.  Thus, in some
sense, what is preventing all Exel-Laca algebras from being
graph algebras is the fact that not all matrices arise as
the edge matrix of a graph.  For example, the finite matrix
$\left( \begin{smallmatrix} 1 & 1 \\ 1 & 0
\end{smallmatrix} \right) $ is not the edge matrix of any
graph.  This is because any such graph would have to have
two edges, $e_1$ and $e_2$, and the relations coming from
the matrix would imply that $r(e_1)=s(e_1)$ and
$r(e_1)=s(e_2)$, but $s(e_1) \neq s(e_2)$.

The way that we will overcome this problem is to allow the
range of each edge to be a set of vertices, rather than just
a single vertex.  For example, if we let $v_1$ and
$v_2$ be two vertices, $e_1$ and $e_2$ be two edges,
and define $s(e_1) = v_1$, $s(e_2) = v_2$, $r(e_1) = \{ v_1,
v_2 \} $, and $r(e_2) = \{ v_1 \}$, then we see that the
``edge matrix" of such an object would be $\left(
\begin{smallmatrix} 1 & 1 \\ 1 & 0
\end{smallmatrix} \right)$ because the edge $e_1$ may be
followed by either $e_1$ or $e_2$, and the edge $e_2$ may
only be followed by $e_1$.  Thus by allowing the edges to
have a set of vertices as their range, we can view the matrix
as an edge matrix.

Recall that a \emph{graph} $G=(G^0,G^1,r,s)$ consists of a
countable set of vertices $G^0$, a countable set of edges
$G^1$, and maps $r,s : G^1 \rightarrow G^0$ identifying the
range and source of each edge.  For a set $X$ let $\Pow(X)$
denote the collection of all subsets of $X$ and let $P(X)$
denote the collection of all nonempty subsets of $X$.

\begin{definition}  An \emph{ultragraph} $\G = (G^0, \G^1, r,
s)$ consists of a countable set of vertices $G^0$, a
countable set of edges $\G^1$, and functions $s : \G^1
\rightarrow G^0$ and $r : \G^1 \rightarrow P(G^0)$. 
\end{definition}

\begin{remark}
Note that an ultragraph is a more general
object than a graph.  A graph may be viewed as a special
type of ultragraph in which
$r(e)$ is a singleton set for each edge $e$.
\end{remark}

\begin{example}  A convenient way to draw ultragraphs is
to first draw the set $G^0$ of vertices, and then for
each edge $e \in \G^1$ draw an arrow labeled $e$ from $s(e)$
to each vertex in $r(e)$.  For instance, the ultragraph given
by
\begin{align*}
G^0 & = \{ v,w,x \} & s(e) & = v & s(f) & = w & s(g) & = x \\
\G^1 & = \{e,f,g \}  & r(e) & = \{ v, w, x \} & r(f) & = \{x
\} & r(g) & = \{v, w \}
\end{align*}
may be drawn as
$$
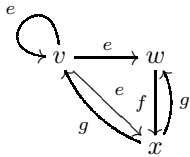
\xymatrix{
v \ar@(u,l)_{e} \ar[r]^{e} \ar[dr]^{e} & w \ar[d]_{f} \\ & x
\ar@/_/[u]_{g} \ar@/^/[ul]^{g} \\ }$$
We then identify any arrows with the same label, thinking of
them as being a single edge.  Thus in the above
example there are only three edges, $e$, $f$,
and $g$, despite the fact that there are six
arrows drawn.
\end{example}

\begin{definition} If $\G$ is an ultragraph, the \emph{edge
matrix} of $\G$ is the $\G^1 \times \G^1$ matrix $A_\G$ given
by $A_\G(e,f) = \begin{cases} 1 & \text{if $s(f) \in r(e)$ }
\\ 0 & \text{otherwise.}  \end{cases}$
\end{definition}

Although not every $\{0,1\}$-matrix is the edge matrix
of a graph, every $\{ 0,1 \}$-matrix is the edge matrix
of an ultragraph.  

\begin{definition}
If $I$ is a countable set and $A$ is an $I
\times I$ matrix with entries in $\{ 0, 1\}$, then we may
form the ultragraph $\G_A := (G_A^0, \G_A^1, r, s)$
defined by $G_A^0 := \{v_i : i \in I \}$, $\G_A^1 :=  I$,
$s(i) = v_i$ for all $i \in I$, and $r(i)=\{v_j : A_\G (i,j)
= 1 \}$.  
\label{edgeLG}
\end{definition}

Note that the edge matrix of $\G_A$ is $A$.  Also note that
for each $v_i \in G_A^0$ there is exactly one edge $i$ with
source $v_i$. 

If $\G$ is an ultragraph, then a vertex $v \in
G^0$ is called a \emph{sink} if $| s^{-1}(v) | = 0$ and an
\emph{infinite emitter} if $| s^{-1} (v) | = \infty$.  We
call a vertex a \emph{singular vertex} if it is
either a sink or an infinite emitter.

\begin{remark}  If $G$ is a graph, then $G$ is said to be
row-finite if it has no infinite emitters.  This terminology
comes from the fact that the edge matrix $A_G$ is
row-finite; that is, the rows of $A_G$ are eventually zero. 
However, this is not the case for ultragraphs:  If $\G$
has no infinite emitters, then it is not necessarily true
that $A_\G$ is row-finite.  In fact, for any matrix $A$ we
see that the ultragraph $\G_A$ will always have no
infinite emitters.  Thus we refrain from using the term
row-finite when speaking of ultragraphs.  Instead
we shall always say that an ultragraph has no infinite
emitters or (in the case that there are also no sinks)
that the ultragraph has no singular vertices.
\end{remark}

For an ultragraph $\G = (G^0, \G^1, r, s)$ we let $\G^0$
denote the smallest subcollection of $\Pow(G^0)$ that
contains $\{v \}$ for all $v \in G^0$, contains $r(e)$ for
all $e \in \G^1$, and is closed under finite intersections
and finite unions. Roughly speaking, the elements of $\{v : v
\in G^0 \} \cup \{r(e) : e \in \G^1 \}$ play the role of
``generalized vertices" and $\G^0$ plays the role of 
``subsets of generalized vertices".

\begin{definition} \label{CK-G-fam} If $\G$ is an ultragraph, a \emph{Cuntz-Krieger $\G$-family} is a collection of
partial isometries $\{ s_e : e \in \G^1 \}$ with mutually
orthogonal ranges and a collection of projections $\{ p_A : A
\in
\G^0
\}$ that satisfy
\begin{enumerate}
\item $p_\emptyset = 0$, $p_A p_B = p_{A \cap B}$, and $p_{A
\cup B} = p_A + p_B - p_{A \cap B}$ for all $A,B \in \G^0$
\item $s_e^*s_e = p_{r(e)}$ for all $e \in \G^1$
\item $s_es_e^* \leq p_{s(e)}$ for all $e \in \G^1$
\item $p_v = \sum_{s(e) = v} s_es_e^*$ whenever $0 < |
s^{-1}(v) | < \infty$.
\end{enumerate}

\noindent When $A$ is a singleton set $\{ v \}$, we write
$p_v$ in place of $p_{ \{ v \} }$.
\end{definition}

For $n \geq 2$ we define $\G^n := \{ \alpha = \alpha_1 \ldots
\alpha_n : \alpha_i \in \G^1 \text{ and } s(\alpha_{i+1}) \in
r(\alpha_i) \}$ and $\G^* := \bigcup_{n=0}^\infty \G^n$.  The
map $r$ extends naturally to $\G^*$, and we say that $\alpha$
has length $|\alpha| = n$ when $\alpha \in \G^n$.  Note that
the paths of length zero are the elements of $\G^0$, and
when $A \in \G^0$ we define $s(A)=r(A)=A$.

\begin{lemma}
If $A \in \G^0$ and $e \in \G^1$, then $$p_A s_e =
\begin{cases} s_e & \text{if $s(e) \in A$} \\ 0 &
\text{otherwise} \end{cases} \hspace{.2in} \text{ and }
\hspace{.2in} s_e^* p_A  = \begin{cases} s_e^* & \text{if
$s(e) \in A$} \\ 0 & \text{otherwise} \end{cases}$$
\label{pslemma}
\end{lemma}

\begin{proof}  We have
$$p_A s_e = p_A s_es_e^*s_e = p_A = p_Ap_{s(e)} s_es_e^*s_e =
p_{A \cap s(e)} s_e = \begin{cases} s_e & \text{if $s(e) \in
A$} \\ 0 & \text{otherwise} \end{cases}$$ and the second claim
follows by taking adjoints.
\end{proof}

For a path $\alpha := \alpha_1 \ldots \alpha_n \in \G^*$ we
define $s_\alpha$ to be $s_{\alpha_1} \ldots s_{\alpha_n}$ if
$|\alpha| \geq 1$ and $p_A$ if $\alpha = A \in \G^0$.

\begin{lemma}
Let $\{ s_e, p_A \}$ be a Cuntz-Krieger $\G$-family, and let
$\beta,\gamma \in \G^*$ with $| \beta|, |\gamma| \geq 1$. 
Then $$s_\beta^*s_\gamma = \begin{cases} s_{\gamma'} &
\text{if $\gamma = \beta \gamma'$, $\gamma' \notin \G^0$} \\
p_{r(\gamma)} & \text{if $\gamma = \beta$} \\ s_{\beta'}^* &
\text{if $\beta = \gamma \beta'$, $\beta' \notin \G^0$}\\ 0 &
\text{otherwise.} \end{cases}$$
\label{sslemma}
\end{lemma}
\begin{proof}
If $e,f \in \G^1$ we have $s_e^*s_f = 0$ unless $e=f$, so
$s_e^*s_\gamma  = \delta_{e,\gamma_1} s_{\gamma_1}^*
s_{\gamma_1} \cdots s_{\gamma_{|\gamma|}} =
\delta_{e,\gamma_1} p_{r(\gamma_1)} s_{\gamma_2} \cdots
s_{\gamma_{|\gamma|}}$.  Because
$s(\gamma_2) \in r(\gamma_1)$, this gives
$s_e^*s_\gamma  = \delta_{e,\gamma_1} s_{\gamma_2} \cdots
s_{\gamma_{|\gamma|}}$.  Repeated calculations of this form
show that $s_\beta s_\gamma = 0$ unless either $\gamma$
extends $\beta$ or $\beta$ extends $\gamma$.  Suppose for the
sake of argument that $\gamma = \beta \gamma'$ extends
$\beta$.  Then calculations as above show that $s_\beta^*
s_\beta s_{\gamma'} = s_{\beta_{|\beta|}}^*
s_{\beta_{|\beta|}} s_{\gamma'} = p_{r(\beta)} s_{\gamma'} =
s_{\gamma'}$.
\end{proof}

\begin{remark}  We see from Lemma~\ref{pslemma} and
Lemma~\ref{sslemma} that any word in $s_e$, $p_A$, and
$s_f^*$ may be written in the form $s_\alpha p_A s_\beta^*$
for some $A \in \G^0$ and some $\alpha, \beta \in \G^*$ with
$r(\alpha) \cap r(\beta) \cap A \neq \emptyset$.   
\label{density} 
\end{remark}

\begin{theorem}
Let $\G$ be an ultragraph.  There exists a
$C^*$-algebra $B$ generated by a universal Cuntz-Krieger
$\G$-family $\{ s_e, p_A \}$.  Furthermore, the $s_e$'s are
nonzero and every $p_A$ with $A \neq \emptyset$ is nonzero.
\end{theorem}

\begin{proof}
We only give an outline here as the argument closely follows
that of \cite[Theorem 2.1]{KPR} and \cite[Theorem 2.1]{aHR}. 
Let $S_\G := \{ (\alpha, A, \beta) : \alpha, \beta \in
\G^*, \ A \in \G^0, \text{ and } r(\alpha) \cap r(\beta) \cap
A \neq \emptyset \}$ and let $k_\G$ be the space of
functions of finite support on $S_\G$.  The set of point
masses $\{ e_\lambda :
\lambda \in S_\G \}$ forms a basis for $k_\G$.  By thinking
of $e_{(\alpha, A, \beta)}$ as $s_\alpha p_A s_\beta^*$ we may
use Lemma~\ref{pslemma}, Lemma~\ref{sslemma}, and the
relation $p_A p_B = p_{A \cap B}$ to define an associative
multiplication and involution on $k_\G$ such that
$k_\G$ is a $*$-algebra.

As a $*$-algebra $k_\G$ is generated by $q_A := e_{(A,A,A)}$
and $t_e := e_{(e, r(e), r(e))}$.  From the way we have
defined multiplication, the elements
$q_A$ have the property that $q_Aq_B = e_{(A,A,A)}e_{(B,B,B)}
= e_{(A \cap B, A \cap B, A \cap B)} = q_{A\cap B}$ and $q_v
\geq t_et_e^*$ for all $e \in \G^1$ with $s(e)=v$.  Let us
mod out by the ideal $J$ generated by the elements $q_v -
\sum_{e : s(e)=v} t_et_e^*$ for all $v$ with $0 <
|s^{-1}(v)| < \infty$ and the elements $q_A+q_B-q_{A \cap B}
- q_{A \cup B}$ for all $A,B \in \G^0$.  Then the images
$r_A$ of $q_A$ and $u_e$ of $t_e$ in $k_\G / J$ form a
Cuntz-Krieger $\G$-family that generates $k_\G / J$.  The
triple $(k_\G, r_A, u_e)$ has the required universal
property, though $k_\G / J$ is not a
$C^*$-algebra.  A standard argument shows that $$\| a \|_0 :=
\sup \ \{ \ \| \pi(a) \| : \pi \text{ is a nondegenerate
$*$-representation of $k_\G / J$} \}$$ is a well-defined,
bounded seminorm on $k_\G / J$.  The completion $B$ of $$(
k_\G / J) / \{ b \in k_\G/J : \| b \|_0 = 0 \}$$ is a
$C^*$-algebra with the same representation theory as $k_\G /
J$.  Thus if $p_A$ and $s_e$ are the images of $r_A$ and
$u_e$ in $B$, then $(B,s_e,p_A)$ has all the required
properties.

Now for each $e \in \G^1$ let $\Hi_e$ be an
infinite-dimensional Hilbert space.  Also for each $v \in
\G^0$ let $\Hi_v := \bigoplus_{s(e)=v} \Hi_e$ if $v$ is not a
sink, and let $\Hi_v$ be an infinite-dimensional Hilbert
space if $v$ is a sink.  Let $\Hi := \bigoplus_{v \in G^0}
\Hi_v$ and for each $e \in \G^1$ let $S_e$ be the partial
isometry with initial space $\bigoplus_{v \in r(e)} \Hi_v$
and final space $\Hi_e$.  Finally, for $A \in \G^0$ define
$P_A$ to be the projection onto $\bigoplus_{v \in A} \Hi_v$,
where this is interpreted as the zero projection when $A =
\emptyset$.  Then $\{ S_e, P_A \}$ is a Cuntz-Krieger
$\G$-family.  By the universal property there exists a
homomorphism $h: B \rightarrow C^*( \{ S_e, P_A \})$.  Since
the $S_e$'s and $P_A$'s are nonzero, it follows that the
$s_e$'s and $p_A$'s are also nonzero. 
\end{proof}
\noindent  The triple $(B,s_e,p_A)$ is unique up to
isomorphism and we write $C^*(\G)$ for $B$.  From
Remark~\ref{density} we see that $C^*(\G) =
\overline{\text{span}} \{ s_\alpha p_A s_\beta^* : \alpha,
\beta \in \G^* \text{ and } A \in \G^0 \}$.  The following
lemma allows us to say slightly more.

\begin{lemma}  If $\G := ( G^0, \G^1,r,s)$ is an ultragraph, then \begin{align*} \G^0 = \{ \bigcap_{e \in X_1} r(e)
\cup \ldots 
\cup \bigcap_{e \in X_n} r(e) \cup F : & \ \text{$X_1,
\ldots, X_n$ are finite subsets of $\G^1$} \\ & \text{ and $F$
is a finite subset of $G^0$} \}.
\end{align*}  
Furthermore, $F$ may be chosen to
be disjoint from $\bigcap_{e \in X_1} r(e) \cup \ldots 
\cup \bigcap_{e \in X_n} r(e)$.
\label{description}
\end{lemma}

\begin{proof}
Recall that $\G^0$ contains  $\{v \}$ for all $v \in G^0$ and
$r(e)$ for all $e \in \G^1$ and is closed under finite
intersections and finite unions.  Hence the right hand side
of the above equation is contained in $\G^0$.  To see the
converse note that the right hand side contains $\{v \}$ for
all $v \in G^0$ and $r(e)$ for all $e \in \G^1$ and is closed
under finite intersections and finite unions.
Furthermore, since $F$ is a finite subset we may always
choose it to be disjoint from $\bigcap_{e \in X_1} r(e) \cup \ldots 
\cup \bigcap_{e \in X_n} r(e)$ simply by discarding any
unwanted vertices in $F$.
\end{proof}

\begin{remark}  This lemma combined with the comment
preceding it shows that \begin{align*} C^*(\G) =
\overline{\text{span}} \{ s_\alpha p_A s_\beta^* : \ & \alpha,
\beta \in \G^* \text{ and either $A = r(e_1) \cap \ldots
\cap r(e_n)$ } \\ & \text{for $e_1 \ldots e_n \in \G^1$ or
$A$ is a finite subset of $\G^0$} \}. 
\end{align*}
Furthermore, one can see that $C^*(\G)$ is generated by
$\{s_e : e \in \G^1 \} \cup \{ p_v : \text{$v$ is a singular
vertex} \}$.
\end{remark}

We conclude with a discussion of the gauge action for
ultragraph algebras.  If $\G$ is an ultragraph and
$\{s_e,p_A \}$ is a Cuntz-Krieger $\G$-family, then for any
$z \in \T$, the family $\{ zs_e, p_A \}$ will be another
Cuntz-Krieger $\G$-family that generates $C^*(\G)$.  Thus
the universal property gives a homomorphism $\gamma_z :
C^*(\G) \rightarrow C^*(\G)$ such that $\gamma_z(s_e) =
zs_e$ and $\gamma_z(p_A) = p_A$.  The homomorphism
$\gamma_{\overline{z}}$ is an inverse for $\gamma_z$, so
$\gamma_z \in \aut C^*(\G)$.  Furthermore, a routine
$\epsilon / 3$ argument shows that $\gamma$ is a strongly
continuous action of $\T$ on $C^*(\G)$.  We call this action
the \emph{gauge action} for $C^*(\G)$.

\section{Viewing graph algebras as ultragraph algebras}
\label{graph-alg-sec}

The construction of $C^*(\G)$ generalizes the
$C^*$-algebra $C^*(G)$ associated to a directed graph $G$ as
described in \cite{KPR} for row-finite graphs and in
\cite{FLR} for arbitrary graphs.  If
$G = (G^0,G^1,r,s)$ is a directed graph, then we may view $G$
as an ultragraph $\G$ in a natural way; that is, let $\G^1
:= G^1$, define $\tilde{r} : \G^1 \rightarrow P(G^0)$ by
$\tilde{r} (e) = \{ r(e) \}$, and then set $\G := (G^0,\G^1,
\tilde{r},s)$.

\begin{proposition}
If $G$ is a graph and $\G$ is the ultragraph associated
to $G$, then $C^*(G)$ is naturally isomorphic to $C^*(\G)$.
\label{graphalg}
\end{proposition}
\begin{proof}  Since $\tilde{r}(e)$ is a singleton set for
all $e \in \G^1$, we see that $\G^0$ equals the collection of
all finite subsets of $G^0$.  If $\{s_e,p_v \}$ is a
Cuntz-Krieger $G$-family \cite{FLR} in $C^*(G)$, then we
may define $p_A := \sum_{v \in A} p_v$ for all $A \in
\G^0$.  (Note that this will be a finite sum.)  Then $\{ s_e,
p_A \}$ is a Cuntz-Krieger $\G$-family.  Conversely,
if $\{t_e, q_A \}$ is a Cuntz-Krieger $\G$-family, then it
restricts to a Cuntz-Krieger $G$-family $\{t_e, q_{\{v\}}
\}$.  The result follows by applying the universal properties.
\end{proof}

\begin{lemma}
If $\G$ is an ultragraph, then $C^*(\G)$ is unital if and
only if $G^0 \in \G^0$, and in this case $1=p_{G^0}$.
\label{unital}
\end{lemma}

\begin{proof}  Let $C^*(\G) = C^*(\{ s_e, p_A \})$.  If $G^0
\in \G^0$, then consider the projection $p_{G^0}$.  For any
$\alpha, \beta \in \G^*$ and $A \in \G^0$ we have $(s_\alpha
p_A s_\beta^*) p_{G^0} = s_\alpha p_A s_\beta^* p_{s(\beta)}
p_{G^0} = s_\alpha p_A s_\beta^* p_{s(\beta)} = s_\alpha p_A
s_\beta^*$.  Similarly, $p_{G^0} (s_\alpha p_A s_\beta^*) =
(s_\alpha p_A s_\beta^*)$.  Since $\{ s_\alpha p_A s_\beta^*
\}$ is dense in $C^*(\G)$, it follows that $p_{G^0}$ is a
unit for $C^*(\G)$.

Conversely, suppose that $C^*(\G)$ is unital.  List the
elements of $\G^1 = \{ e_1, e_2, \ldots \}$ and $G^0
= \{v_1, v_2, \ldots \}$.  Note that these sets are either
finite or countably infinite.  For $n
\geq 1$ define $A_n := \bigcup_{i=1}^n \{v_i\} \cup
\bigcup_{i=1}^n r(e_i)$.  Then $A_n \in \G^0$ for all $n$ and
$A_1 \subseteq A_2 \subseteq \ldots$.  Also $\{p_{A_i} \}$ is
an approximate unit since for any $s_\alpha p_A s_\beta^*$ we
may choose $n$ large enough so that $p_{A_n}$ acts as the
identity on $s_\alpha p_A s_\beta^*$.  Now if $A_n
\subsetneq A_m$ for some $m >n$, then there exists $v \in
A_m \backslash A_n$ and
$(p_{A_m}-p_{A_n} )p_v = p_{A_m \cap \{v\}}-p_{A_n \cap \{ v
\}} = p_v$ so
$p_{A_m}-p_{A_n} \geq p_v$ and $\| p_{A_m}-p_{A_n} \| =1$. 
Since $C^*(\G)$ is unital, we must have $p_{A_n} \rightarrow
1$ in norm.  But the only way that this could happen is if
$p_{A_n}$ is eventually constant.  Thus $p_{A_k} = 1$ for
some $k$.  Furthermore, $A_k$ must be all of $G^0$ for if
there were $v \in G^0 \backslash A_k$, then $p_{A_k} p_v = 0$
contradicting the fact that $p_{A_k}=1$.  Hence $A_k = G^0$
and since $A_k \in \G^0$ we are done.
\end{proof}

\begin{remark}  Note that when a graph $G$ is viewed as an
ultragraph as in Proposition~\ref{graphalg}, then the above
lemma produces the familiar result that $C^*(G)$ is unital
if and only if $G$ has a finite number of vertices.
\end{remark}

Recall that a loop in a graph $G$ is a path $\alpha \in G^*$
with $| \alpha | \geq 1$ and $s(\alpha) = r(\alpha)$.  If
$\alpha = \alpha_1 \ldots \alpha_n$ is a loop, then an exit
for $\alpha$ is defined to be an edge $e \in G^1$ with $s(e)
= s(\alpha_i)$ for some $1 \leq i \leq n$ but $e \neq
\alpha_i$.  A graph $G$ is said to satisfy Condition~(L) if
all loops have exits.  Roughly speaking, an exit for a loop
$\alpha:= \alpha_1
\ldots \alpha_n$ is something that allows you to ``get out"
of the loop $\alpha$; that is, it allows you to follow a
path other than $\alpha_1 \ldots \alpha_n \alpha_1 \ldots
\alpha_n \alpha_1 \ldots$. 

\begin{definition}
If $\G$ is an ultragraph, then a \emph{loop} is a path
$\alpha \in \G^*$ with $| \alpha | \geq 1$ and $s(\alpha)
\in r(\alpha)$.  An \emph{exit} for a loop is either of the
following:
\begin{enumerate}
\item an edge $e \in
\G^1$ such that there exists an $i$ for which $s(e) \in
r(\alpha_{i})$ but $e \neq \alpha_{i+1}$
\item a sink $w$ such that $w \in r(\alpha_i)$ for some $i$.
\end{enumerate}
\end{definition}

\noindent We now extend Condition~(L) to ultragraphs.

\text{ }

\noindent \textbf{Condition~(L):}  Every loop in $\G$ has an
exit; that is, for any loop $\alpha := \alpha_1 \ldots
\alpha_n$ there is either an edge $e \in \G^1$ such that
$s(e) \in r(\alpha_{i})$ and $e \neq
\alpha_{i+1}$ for some $i$, or there is a sink $w$ with $w
\in r(\alpha_i)$ for some $i$.

\text{ }

\noindent Note that if $\alpha := \alpha_1 \ldots \alpha_n$
is a loop in $\G$ without an exit, then for all $i$ we must
have $r(\alpha_i) = \{ s(\alpha_{i+1}) \}$ and $s^{-1}
(s(\alpha_i)) = \{ \alpha_i \}$.

\section{Viewing Exel-Laca algebras as ultragraph
algebras}
\label{EL-sec}

In this section we shall see that the $C^*$-algebras of
ultragraphs with no singular vertices are precisely the
Exel-Laca algebras.

\begin{definition}[Exel-Laca]  Let $I$ be any set and let $A =
\{A(i,j)_{i,j \in I} \}$ be a $\{0,1\}$-matrix
over $I$ with no identically zero rows.  The Exel-Laca
algebra $\OA$ is the universal
$C^*$-algebra generated by partial isometries $\{ s_i : i \in
I \}$ with commuting initial projections and mutually
orthogonal range projections satisfying $s_i^* s_i s_j s_j^*
= A(i,j) s_js_j^*$ and 
\begin{equation} 
\label{conditionfour}
\prod_{x \in X} s_x^*s_x \prod_{y \in Y} (1-s_y^*s_y) =
\sum_{j \in I} A(X,Y,j) s_js_j^* 
\end{equation}
whenever $X$ and $Y$ are finite subsets of $I$ such that the
function $$j \in I \mapsto A(X,Y,j) := \prod_{x \in X} A(x,j)
\prod_{y \in Y} (1-A(y,j))$$ is finitely supported.
\end{definition}

Although there is reference to a unit in
(\ref{conditionfour}), this relation applies to algebras that
are not necessarily unital, with the convention that if a $1$
still appears after expanding the product in
(\ref{conditionfour}), then the relation implicitly states
that
$\OA$ is unital.  It is also important to realize that the
relation (\ref{conditionfour}) also applies when the function
$j \mapsto A(X,Y,j)$ is identically zero.  This particular
instance of (\ref{conditionfour}) is interesting in itself so
we emphasize it by stating the associated relation
separately:
\begin{equation}
\label{speccondfour}
\prod_{x \in X} s_x^*s_x \prod_{y \in Y} (1-s_y^*s_y) = 0
\end{equation}
whenever $X$ and $Y$ are finite subsets of $I$ such that
$A(X,Y,j)=0$ for every $j \in I$.

\begin{lemma}
Let $\G$ be an ultragraph.  If $A \subseteq G^0$ is a finite
set, then $$p_A = \sum_{v \in A} p_v.$$
\label{sum}
\end{lemma}
\begin{proof} Simply use the fact that $A$ is the disjoint
union of its singleton sets.
\end{proof}

\begin{lemma}
Let $\G$ be an ultragraph.  If $Y \subseteq \G^1$ is a
finite set, then for any $A \in \G^0$ $$\prod_{y \in Y} (p_A
-p_A p_{r(y)}) = p_A -p_A p_B \hspace{.2in} \text{ where $B =
\bigcup_{y \in Y} r(y)$.}$$
\label{prod}
\end{lemma}
\begin{proof}
We shall induct on the number of elements in $Y$.  If $| Y |
= 1$, then the claim holds trivially.  Assume the claim is
true for sets containing $n-1$ elements.  Let $| Y | = n$ and
choose $e \in Y$.  If we let $B' := \bigcup_{y \in Y
\backslash \{ e \}} r(y)$, then
\begin{align*}\prod_{y
\in Y} (p_A -p_A p_{r(y)}) & = \left( \prod_{y \in Y
\backslash \{ e \} } (p_A -p_A p_{r(y)}) \right) (p_A - p_A
p_{r(e)}) \\ & = (p_A - p_Ap_{B'})(p_A - p_A p_{r(e)}) \\ & =
p_A - p_A p_{B'} - p_A p_{r(e)} + p_A p_{B' \cap r(e)} \\ & =
p_A - p_A (p_{B'} + p_{r(e)} - p_{B' \cap r(e)}) \\ & = p_A
- p_A p_{B' \cup r(e)} \\ & = p_A-p_Ap_B
\end{align*}
where $B := \bigcup_{y \in Y} r(y)$.
\end{proof}

\begin{proposition}
Let $\G$ be an ultragraph with no sinks, and
let $A_\G$ be the edge matrix of $\G$.  If $\{s_e, p_A \}$ is
a Cuntz-Krieger $\G$-family, then $\{ s_e : e \in \G^1 \}$ is
a collection of partial isometries satisfying the relations
defining $\mathcal{O}_{A_\G}$.
\label{edgematrix}
\end{proposition}
\begin{proof}
By definition the $s_e$'s have mutually orthogonal
range projections.  Furthermore, $s_e^*s_es_f^*s_f =
p_{r(e)} p_{r(f)} = p_{r(e) \cap r(f)} = p_{r(f)} p_{r(e)} =
s_f^*s_fs_e^*s_e$ so the initial projections commute. 
Furthermore, if $e,f \in \G^1$, then by Lemma~\ref{pslemma}
$$s_e^*s_es_fs_f^* = p_{r(e)} s_fs_f^* = \begin{cases}
 s_fs_f^* & \text{if $s(f) \in r(e)$} \\ 0 &
\text{otherwise} \end{cases} = A_\G (e,f) s_fs_f^*.$$

Now we shall show that condition (\ref{conditionfour}) of
Exel-Laca algebras holds.  Suppose that $X$ and $Y$ are
finite subsets of $\G^1$ such that the function $j
\mapsto A_\G(X,Y,j)$ has finite (or empty) support.  We
divide the proof of (\ref{conditionfour}) into two cases.

\noindent \textsc{Case 1:}  $X = \emptyset$.

We claim that $G^0 \in \G^0$.  The support of $j
\mapsto \prod_{y \in Y} (1-A_{\G}(y,j))$, which is finite (or
empty) by assumption, is given by $$F := \{ j \in \G^1 :
A_{\G}(y,j) = 0 \text{ for all $y \in Y$} \} = \{ j \in \G^1
: s(j) \notin \bigcup_{y \in Y} r(y) \}.$$  Because there are
no sinks, the source map $s$ is surjective and $G^0 = s(F)
\sqcup \bigcup_{y \in Y} r(y)$.  Since $F$ and $Y$ are finite
this implies that $G^0 \in \G^0$ and by Lemma~\ref{unital}
$C^*(\G)$ is unital with $1=p_{G^0}$.  Furthermore, since $G^0
= s(F) \sqcup \bigcup_{y \in Y} r(y)$  we have $$1 =
p_{s(F)} + p_B \hspace{.2in} \text{ where } B := \bigcup_{y
\in Y} r(y)$$ and using Lemma~\ref{sum} $$1-p_B = \sum_{v \in
s(F)} p_v.$$  Now applying Lemma~\ref{prod} with $p_A = 1$
gives $$\prod_{y \in Y} ( 1 - p_{r(y)}) = \sum_{v \in s(F)}
p_v$$ or $$\prod_{y \in Y} ( 1 - s_y^*s_y) = \sum_{v \in s(F)}
p_v.$$  If $F$ is empty the right hand side vanishes and
(\ref{conditionfour}) holds.  If $F$ is nonempty, then $s(F)
\neq \emptyset$ and for each $v \in s(F)$ we have $\{ j :
s(j) = v \} \subseteq \{ j : j \notin \bigcup_{y \in Y} r(y)
\} = F$.  Hence $0 < | \{ j : s(j) = v \} | < \infty$ and we
may use the definition of a Cuntz-Krieger $\G$-family to
write $p_v = \sum_{ \{ j : s(j) = v \} } s_js_j^*$.  Summing
over all vertices in $s(F)$ and substituting above gives
$$\prod_{y \in Y} (1 - s_y^*s_y) = \sum_{j \in F} s_js_j^*$$
which is (\ref{conditionfour}).

\noindent \textsc{Case 2:} $X \neq \emptyset$.

Once again let $F = \{ j \in \G^1 : s(j) \in r(x) \text{
for all } x \in X \text{ and } s(j) \notin r(y) \text{
for all } y \in Y \}$ be the support of $j \mapsto
A(X,Y,j)$.
\begin{itemize}
\item  \textsc{Subcase 2}a:  Assume $F$ is empty.  Let $x_0
\in X$.  Since there are no sinks, for each $v \in r(x_0)$ 
there exists $j_v \in \G^1$ such that $s(j_v) = v$.  Because
$F$ is empty, either $v = s(j_v) \notin r(x_v)$ for some $x_v
\in X$ or $v = s(j_v) \in r(y_v)$ for some $y_v \in Y$. 
Thus $$r(x_0) \cap \bigcap_{x_v} r(x_v) \subseteq
\bigcup_{y_v} r(y_v)$$ and the left hand side of
(\ref{speccondfour}) contains $$\left( \prod_{x_v}
s_{x_v}^*s_{x_v} \right) \left( s_{x_0}^* s_{x_0} \right)
\prod_{y_v} (1-s_{y_v}^*s_{y_v})  = \prod_{x_v} p_{r(x_v)} 
\prod_{y_v} \left( p_{r(x_0)}-p_{r(x_0)}p_{r(y_v)} \right).$$ 
If we let $A := \bigcap_{x_v} r(x_v)$ and $B := \bigcup_{y_v}
r(y_v)$, then Lemma~\ref{prod} shows that 
 
\begin{align*} \left( \prod_{x_v} s_{x_v}^*s_{x_v} \right)
\left( s_{x_0}^* s_{x_0} \right)
\prod_{y_v} (1-s_{y_v}^*s_{y_v}) & = p_A (p_{r(x_0)}
- p_{r(x_0)} p_B) \\
& = p_{r(x_0) \cap A} - p_{r(x_0) \cap A \cap B} \\
& = p_{r(x_0) \cap A} - p_{r(x_0) \cap A} \\
& = 0
\end{align*}
so (\ref{speccondfour}) holds.

\item \textsc{Subcase 2}b: Suppose $F$ is nonempty.  Since
$F =
\{ j : s(j) \in \bigcap_{x \in X} r(x) \cap \left( \bigcup_{y
\in Y} r(y) \right)^C \}$ is nonempty and finite and since
$\G$ has no sinks, it follows that $\bigcap_{x \in X} r(x)
\cap
\left( \bigcup_{y \in Y} r(y) \right)^C$ is nonempty and
finite.  Thus  $\bigcap_{x \in X} r(x) \cap
\left( \bigcup_{y \in Y} r(y) \right)^C \in \G^0$.  For
convenience of notation, let $A := \bigcap_{x \in X} r(x)$ and
$B := \bigcup_{y \in Y} r(y)$.  Then
$$\sum_{j \in F} s_js_j^* = \sum_{ \{j : s(j) \in A \cap B^C
\} } s_js_j^*.$$
Furthermore, since $F = \{ j : s(j) \in A \cap B^C \}$ is
finite, we know that any vertex in $A \cap B^C$ can emit only
finitely many vertices.  Thus $p_v = \sum_{s(e) = v}
s_es_e^*$ for all $v \in A \cap B^C$.  Applying
this to the above equation and using Lemma~\ref{sum} gives
$$\sum_{j \in F} s_js_j^* = \sum_{ v \in A \cap B^C}
p_v = p_{A \cap B^C}.$$  Now we also have $$p_{A \cap B^C} +
p_{A \cap B} = p_A - p_\emptyset = p_A$$ and combining this
with Lemma~\ref{prod} gives
\begin{align*}
\sum_{j \in F} s_js_j^* & = p_A - p_{A \cap B} \\ &= 
\prod_{y \in Y} (p_A - p_A  p_{r(y)}) \\ & = \prod_{x \in X}
p_{r(x)} \prod_{y \in Y} (1- p_{r(y)}) \\ &= \prod_{x \in X}
s_x^*s_x \prod_{y \in Y} (1- s_y^*s_y)
\end{align*}
so (\ref{conditionfour}) holds and we are done.
\end{itemize}
\end{proof}

\begin{theorem}
Let $\G$ be an ultragraph with no singular vertices (i.e.~no
sinks and no infinite emitters).  If $A_\G$ is the edge
matrix of $\G$, then $\mathcal{O}_{A_\G}$ is canonically
isomorphic to $C^*(\G)$.
\label{ELareLG} 
\end{theorem}

\begin{proof}  By Proposition~\ref{edgematrix} and the
universal property of $\mathcal{O}_{A_\G}$, there exists a
homomorphism $\phi : \mathcal{O}_{A_\G} \rightarrow C^*(\G)$
with the property that $\phi (S_e) = s_e$.  Since this
homomorphism is equivariant for the gauge actions, it
follows from the Gauge-Invariant Uniqueness Theorem for
Exel-Laca algebras
\cite[Theorem 2.7]{RS} that $\phi$ is injective. 
Furthermore, since $\G$ has no singular vertices, the
$s_e$'s generate $C^*(\G)$.  Thus $\phi$ is also surjective.
\end{proof}

\begin{remark}
Since any matrix $A$ with entries in $\{0,1\}$ is the edge
matrix of the ultragraph $\G_A$, the above shows that the
$C^*$-algebras of ultragraphs with no singular vertices
are precisely the Exel-Laca algebras.
\end{remark}

\begin{corollary}
Let $\G$ be an ultragraph with no singular vertices.  If
$A$ is the edge matrix of $\G$, then $C^*(\G)$
is canonically isomorphic to $C^*(\G_A)$.
\end{corollary}
\begin{proof} From Theorem~\ref{ELareLG} we have $C^*(\G)
\cong \OA \cong C^*(\G_A)$.
\end{proof}

\begin{remark}
This corollary shows that if one wishes to
study ultragraphs with no singular vertices, then there
is no loss of generality in considering only ultragraphs
of the form
$\G_A$ for a matrix $A$ with entries in $\{0,1\}$.  This is
somewhat surprising since any $\G_A$ has the property that $|
s^{-1}(v) | = 1$ for all $v \in G^0$.
\end{remark}

In conclusion, we have seen that if $A$ is a
$\{0,1\}$-matrix with no zero rows, then $\OA \cong
C^*(\G_A)$, and hence we may view any Exel-Laca algebra as an ultragraph algebra.  While this has the advantage that
one no longer needs to deal with the complicated
Condition~\ref{conditionfour} of Exel-Laca algebras, one
must now deal with the collection
$\G_A^0$.  At first it may seem as though we are simply
trading one set of troubles for another.  However, we
contend that the ultragraph approach is often easier. 
Despite the somewhat complicated definition of $\G^0$, we
shall see in the next sections that techniques much like
those used for graph algebras can be applied to ultragraph algebras.

\section{Uniqueness Theorems for Labeled Graph Algebras}
\label{uniqueness-sec}

In this section we prove versions of the Gauge-Invariant
Uniqueness Theorem and the Cuntz-Krieger Uniqueness Theorem
for $C^*$-algebras of ultragraphs with no singular
vertices.  In \S\ref{sing-sec} we extend these results to
allow for singular vertices.  We obtain our results in this
section by approximating ultragraph algebras with
$C^*$-algebras of finite graphs as in \cite{RS}.  We mention
that Raeburn and Szyma\'nski gave one method for
approximating $C^*$-algebras of infinite graphs by those of
finite graphs \cite[Definition 1.1]{RS} and another for
approximating Exel-Laca algebras by $C^*$-algebras of finite
graphs \cite[Definition 2.1]{RS}.  We shall show that we are
able to approximate $C^*$-algebras of ultragraphs by a
method much like the one used for Exel-Laca algebras.  We
also remark that although our computations in
Proposition~\ref{approx} are done for ultragraphs rather
than matrices, they are similar to those in
\cite[Proposition 2.2]{RS}.

Let $\G$ be an ultragraph.   For any finite subsets
$X,Y \subseteq \G^1$ define $V(X,Y) := \bigcap_{x \in X} r(x)
\cap ( \bigcup_{y \in Y} r(y))^c$ and 
$E(X,Y) := \{ e \in \G^1 : s(e) \in V(X,Y) \}$.

\begin{definition} For any finite subset $F \subseteq \G^1$
define the graph $G_F$ by
\begin{align*} G_F^0 & := F \cup \{ X : \emptyset
\neq X \subseteq F \text{ satisfies } E(X,F \backslash X)
\nsubseteq F \}, \text{ and} \\ G_F^1 &:= \{ (e,f) \in F
\times F : s(f) \in r(e) \} \cup \{ (e,X) : e \in X \}; 
\end{align*}
with
\begin{align*}
s((e,f)) & = e & s((e,v)) & = e & s((e,X)) & = e \\
r((e,f)) & = f & r((e,v)) & = v & r((e,X)) & = X.
\end{align*}
\end{definition}

\begin{lemma}  If $P_1, \ldots, P_n$ are commuting
projections, then $$ 1 = \sum_{Y \subseteq \{ 1, \ldots n \}
} \Big( \prod_{i \in Y} P_i \Big) \Big( \prod_{i \notin Y}
(1-P_i) \Big).$$ \label{sumpro}
\end{lemma}

\begin{proof} Induct on $n$: multiply the formula
for $n=k$ by $P_{k+1} + (1-P_{k+1})$. \end{proof}

\begin{proposition}  Let $\G$ be an ultragraph with no
sinks, $\{ s_e, p_A \}$ be a Cuntz-Krieger
$\G$-family, and $F$ a finite subset of $\G^1$.  Define $A_X
:= \bigcap_{x \in X} r(x)$ and
$B_X := \bigcup_{y \in F \backslash X} r(y)$.  Then
\begin{align*}  
Q_e  & := s_es_e^*  &  Q_X & := p_{A_X} (1- p_{B_X}) \Big(
1- \sum_{f \in F} s_fs_f^* \Big) \\  
T_{(e,f)}  & := s_eQ_f   &  T_{(e,X)} & := s_eQ_X  
\end{align*}
forms a Cuntz-Krieger $G_F$-family that generates $C^*(
\{ s_e : e \in  F \} )$.  If every $s_e$ is
nonzero, then every $Q$ is nonzero.
\label{approx}
\end{proposition}

\begin{proof}
We shall first show that this is a Cuntz-Krieger
$G_F$-family.  The projections $Q_e$ are mutually
orthogonal because the $s_e$'s have mutually orthogonal
ranges, and are orthogonal to the $Q_X$'s because of the
factor $1- \sum_{f \in F} s_fs_f^*$.  To see that the $Q_X$'s
are mutually orthogonal suppose that $X \neq Y$.  Then,
without loss of generality, we may assume that there exists
$x \in X \backslash Y$, and because $X \subseteq F$ we see
that $A_X \subseteq r(x)$ and $r(x) \subseteq B_Y$.  Thus
$p_{A_X} (1-p_{B_Y}) = p_{A_X} p_{r(x)}(1- p_{B_Y}) =
p_{A_X} (p_{r(x)} - p_{r(x)} p_{B_Y})  = p_{A_X} (p_{r(x)} -
p_{r(x)} ) = 0$ and $Q_X$ is orthogonal to $Q_Y$.  

\noindent Furthermore,
\begin{align*}
T_{(e,f)}^*T_{(e,f)} & = Q_f^*s_e^*s_eQ_f = s_fs_f^*p_{r(e)}
s_fs_f^* = s_fs_f^*p_{s(f)} p _{r(e)} s_fs_f^* \\
&  = s_fs_f^*p_{s(f)} s_fs_f^* = s_fs_f^* = Q_f,
\end{align*}
and since $Q_X \leq s_e^*s_e$ whenever $e \in X$, we have 
$$ T_{(e,X)}^*T_{(e,X)}  = Q_X^* s_e^* s_e Q_X = Q_X $$
so the first Cuntz-Krieger relation holds.

Note that the elements $X$ in $G_F^0$ are
all sinks in $G_F$.  Thus to check the second
Cuntz-Krieger relation, we need only consider edges whose
source is some $e \in F$.

If $X$ is a subset of $G^1$ and $E(X,F\backslash
X) \subseteq F$, then $V(X,F \backslash X) := \bigcap_{x \in
X} r(x) \cap ( \bigcup_{y \in (F \backslash X)} r(y))^c$ is
a finite set and $V(X,Y) \in \G^0$.  Furthermore, since
$E(X,F \backslash X) \subseteq F$ we see that each vertex in
$V(X,F \backslash X)$ is the source of finitely many edges. 
Hence $$p_{V(X,F \backslash X)} = \sum_{v \in V(X,F
\backslash X)} p_v = \sum_{e \in E(X,F \backslash X)}
s_es_e^* \leq \sum_{f \in F} s_fs_f^*$$ and
$$p_{A_X} (1-p_{B_X}) \Big (1- \sum_{f \in F} s_fs_f^* \Big)
= p_{V(X,Y)} \Big( 1- \sum_{f \in F} s_fs_f^* \Big) =  0.$$ 
Thus if we fix $e \in G^1$, we have
\begin{align}
\sum_{ \{ X : e \in X \} } Q_X 
& =  \sum_{ \{ X : e \in X \} } p_{A_X} (1-p_{B_X} ) \Big( 1-
\sum_{f \in F} s_fs_f^* \Big)   \notag \\   
& = p_{r(e)} \Big( \sum_{ Y \subseteq F
\backslash \{e \} } p_{A_Y} (1-p_{B_{Y \cup \{e \}}} )
\Big( 1- \sum_{f \in F} s_fs_f^* \Big) \Big)  \notag \\
& \quad \quad \text{and by Lemma~\ref{prod}} \notag \\
& = p_{r(e)}  \sum_{ Y
\subseteq F \backslash \{e \} } \Big( \prod_{x \in Y}
p_{r(x)} \Big)  \Big( \prod_{y \in (F \backslash \{e \} )
\backslash Y} (1-p_{r(y)})  \Big)  \Big( 1- \sum_{f
\in F} s_fs_f^* \Big) \notag \\
& \quad \quad \text{and by Lemma~\ref{sumpro}}  \notag \\
& = p_{r(e)} \Big( 1- \sum_{f \in F} s_fs_f^* \Big) \notag \\
& = s_e^*s_e \Big( 1 - \sum_{ \{f \in F : s(f)
\in r(e) \} } s_fs_f^* \Big).  \label{Esum}
\end{align}

\noindent Now we have 
$$\sum_{s(f) \in r(e)} T_{(e,f)}T_{(e,f)}^* + \sum_{ \{ X :
e \in X \} } T_{(e,X)}T_{(e,X)}^* = \sum_{s(f) \in r(e)}
s_es_fs_f^*s_e^* + \sum_{ \{ X : e \in X \} } s_eQ_Xs_e^*$$
which equals $s_es_e^* = Q_e$ by (\ref{Esum}).  Thus the
$T$'s and $Q$'s form a Cuntz-Krieger $G_F$-family.

Equation (\ref{Esum}) also implies that we can recover $s_e$
as $$s_e = \sum_{
\{f \in F : s(f) \in r(e) \} } T_{(e,f)} + \sum_{ \{ X
: e \in X \} } T_{(e,X)} = s_e \Big( \sum_{ \{ X
: e \in X \} } Q_X + \sum_{ \{f \in F : s(f) \in r(e) \} }
s_fs_f^* \Big)$$ so the operators
$T_e$ and $Q_v$ generate $C^*(\{ s_e : e \in F \})$.  For
the last comment note that $E(X,F \backslash X ) \nsubseteq
F$ implies $Q_X \geq s_fs_f^*$ for some $f \notin F$, and
hence $Q_X \neq 0$.
\end{proof}

\begin{corollary}
Let $\G$ be an ultragraph with no sinks and let $F
\subseteq \G^1$ be a finite set of edges.  Then $C^*(G_F)$
is canonically isomorphic to the $C^*$-subalgebra of
$C^*(\G)$ generated by $\{ s_e : e \in F \}$.
\label{subalgebra}
\end{corollary}
\begin{proof}
Applying the proposition to the canonical family
$\{s_e,p_A \}$ of $C^*(\G)$ gives a Cuntz-Krieger
$G_F$-family $\{T_e,Q_v \}$ that generates $C^*(\{ s_e : e
\in F \})$.
Thus we have a homomorphism $\phi:C^*(G_F) \rightarrow
C^*(\G)$ whose image is $C^*( \{ s_e : e \in F \})$. If
$\alpha$ is the gauge action on $C^*(G_F)$ and
$\gamma$ is the gauge action on $C^*(\G)$, then we see
that $\phi \circ \alpha_z = \gamma_z \circ \phi$ for all
$z \in \T$.  Since each projection $Q_v$ is nonzero, it
follows from the Gauge-Invariant Uniqueness Theorem for
graph algebras \cite[Theorem 2.1]{BPRS} that $\phi$ is
injective.
\end{proof}

The following is an analogue of the Gauge-Invariant Uniqueness Theorem for $C^*$-algebras of ultragraphs with
no singular vertices.  We shall extend this result to all
ultragraph algebras in \S\ref{sing-sec}.

\begin{proposition}
Let $\G$ be an ultragraph with no singular vertices (i.e.
no sinks or infinite emitters).  Also let
$\{ S_e,P_A \}$ be a Cuntz-Krieger $\G$-family on Hilbert
space and let $\pi$ be the representation of $C^*(\G)$ such
that $\pi (s_e) = S_e$ and $\pi(p_A) = P_A$.  Suppose that
each $P_A$ is nonzero for every nonempty $A$, and that there
is a strongly continuous action $\beta$ of $\T$ on
$C^*(S_e,P_A)$ such that $\beta_z \circ \pi = \pi \circ
\gamma_z$ for all $z \in \T$.  Then $\pi$ is faithful.
\label{gauge-invariant-nosing}
\end{proposition}

\begin{proof}
Let $F$ be a finite subset of $\G^1$.  Then $C^*(\{ s_e : e
\in F \} )$ is isomorphic to the graph algebra $C^*(G_F)$ by
Corollary \ref{subalgebra}, and this isomorphism is
equivariant for the gauge actions.  Furthermore, the
projections in $\B ( \Hi)$ corresponding to the vertices of
$G_F$ are all nonzero: $Q_e$ because $T_e$ is, and $Q_X$
because of the existence of an $f$ such that $T_fT_f^* \leq
Q_X$.  Applying the Gauge-Invariant Uniqueness Theorem for graph algebras \cite[Theorem
2.1]{BPRS} to the corresponding representation of $C^*(G_F)$
shows that $\pi$ is faithful on $C^*(\{s_e : e \in F \})$,
and hence is isometric there.  Thus $\pi$ is isometric on
the subalgebra generated by $\{ s_e : e \in \G^1 \}$.  Since
$\G$ has no singular vertices, this subalgebra is dense in
$C^*(\G)$.  Hence $\pi$ is isometric on all of $C^*(\G)$.
\end{proof}

We can also use this method of approximating ultragraph
algebras by graph algebras to prove a version of the
Cuntz-Krieger Uniqueness Theorem.  Again, we shall extend
this result to all ultragraph algebras in
\S\ref{sing-sec}.

\begin{lemma}
Let $\G$ be ultragraph with no sinks and let $F \subseteq
\G^1$ be a finite set.  If $L = x_1 \ldots x_n$ is a loop in
$G_F$, then there exists a loop $L' = e_1 \ldots e_n$ in
$\G$ such that $\{e_i \}_{i=1}^n \subseteq F$, $x_i =
(e_i,e_{i+1})$ for $i=1,\ldots,n-1$ and $x_n = (e_n,e_1)$. 
Furthermore,
$L$ has an exit if and only if $L'$ does.
\label{loopexit}
\end{lemma}

\begin{proof}  If $L$ has an exit, then either there is an
edge of the form $(e_i,f)$ with $f \neq e_{i+1}$, or there
is an edge of the form $(e_i,X)$.  In the first case, $s(f)
\in r(e_i)$ and $f \neq e_{i+1}$ so $f$ is an exit for
$L'$.  In the second case $e_i \in X$ and $E(X,F \backslash
X) \nsubseteq F$.  Hence there exists $g \in \G^1$ such that
$s(g) \in \G^1 \backslash F$ for which $s(g) \in V(X,F
\backslash F)$.  But then $s(g) \in r(e_i)$ and $g \neq
e_{i+1}$ so $g$ is an exit for $L'$.

Conversely, suppose $L'$ has an exit.  Since $\G$ has no
sinks, there exists an edge $f \in \G^1$ with $s(f) \in
r(e_i)$ and $f \neq e_{i+1}$.  If $f \in F$, then $(e_i,f)$
is an exit for $L$.  If $f \notin F$, let $X := \{ e \in F :
s(f) \in r(e) \}$.  Then $f \in E(X,F \backslash F)$, so
$\emptyset \neq E(X,F \backslash X) \nsubseteq F$.  Since
$e_i \in X$ we see that $(e_i,X)$ is an exit for $L$.
\end{proof}

\begin{proposition}
Suppose that $\G$ is an ultragraph with no singular
vertices (i.e. no sinks or infinite emitters), and that
$\G$ satisfies condition~(L).  If
$\{S_e,P_A \}$ and $\{T_e, Q_A \}$ are two Cuntz-Krieger
$\G$-families in which all the projections $P_A$ and $Q_A$
are nonzero for nonempty $A$, then there is an isomorphism
$\phi$ of $C^*(S_e,P_A)$ onto $C^*(T_e,Q_A)$ such that
$\phi(S_e)=T_e$ and $\phi(P_A) = Q_A$.
\label{uniqueness-nosing}
\end{proposition}

\begin{proof}
We shall prove the theorem by showing that the
representations $\pi_{S,P}$ and
$\pi_{T,Q}$ of $C^*(\G)$ are faithful, and then $\phi :=
\pi_{T,Q} \circ \pi_{S,P}^{-1}$ is the required isomorphism.

Write $\G^1 = \bigcup_{n=1}^\infty F_n$ as the increasing
union of finite subsets $F_n$, and let $B_n$ be the
$C^*$-subalgebra of $C^*(\G)$ generated by $\{ s_e : e \in
F_n \}$.  By Lemma~\ref{subalgebra} there are isomorphisms
$\phi_n : C^*(G_{F_n}) \rightarrow B_n$ that respect the
generators.  Since all the loops in $F_n$ have exits by
Lemma~\ref{loopexit}, the Cuntz-Krieger Uniqueness Theorem
for graph algebras \cite[Theorem 3.1]{BPRS} implies that
$\pi_{S,P} \circ \phi_n$ is an isomorphism, and hence is
isometric.  Thus $\pi_{S,P}$ is isometric on the
$*$-subalgebra $\bigcup_{n} B_n$ of $C^*(\G)$.  But since
$\G$ has no sinks or infinite emitters, $C^*(\G)$ is
generated by the $s_e$'s and thus $\bigcup_{n} B_n$ is a
dense $*$-subalgebra of
$C^*(\G)$.  Hence $\pi_{S,P}$ is isometric on all of
$C^*(\G)$, and in particular, it is an isomorphism.
\end{proof}

Let $I$ be a countable (or finite) set, and let $A$
be an $I \times I$ matrix with entries in $\{ 0 , 1\}$ and
no zero rows.  In \cite{EL} Exel and Laca associated a graph
$\Gr(A)$ to $A$ whose vertex matrix is equal to $A$. 
Specifically, we define the vertices of $\Gr (A)$ to be $I$,
and for each pair of vertices $i,j \in I$ we define there to
be $A(i,j)$ edges from $i$ to $j$.

In \cite[\S 13]{EL} Exel and Laca proved a uniqueness theorem
for $\OA$ when $\Gr(A)$ satisfies Condition~(L) (or in their
terminology, when $\Gr(A)$ has no terminal circuits).  The
following shows that their uniqueness theorem is equivalent
to the one we proved in Proposition~\ref{uniqueness-nosing}. 

\begin{lemma}
Let $\G$ be an ultragraph with no sinks and with edge
matrix $A$.  Then $\Gr(A)$ satisfies Condition~(L) if and
only if $\G$ satisfies Condition~(L). 
\label{condL}
\end{lemma}
\begin{proof}
Suppose that $\G$ satisfies Condition~(L).  Let $\alpha :=
\alpha_1 \ldots \alpha_n$ be a loop in $\Gr(A)$ with $a_i :=
s(\alpha_i)$ for $1 \leq i \leq n$.  Then $a_1 \ldots a_n$
is a loop in $\G$.  Let $b$ be an exit for this loop, and
without loss of generality assume that $s(b) \in r(a_1)$ and
$b \neq a_2$.  Since $A(a_1,b)=1$, there exists an edge $f$
in $\Gr(A)$ from $a_1$ to $b$.  Since $b \neq a_2$ we know
that $f \neq \alpha_1$ and hence $f$ is an exit for $\alpha$.

Conversely, suppose that $\Gr(A)$ satisfies Condition~(L). 
Let $a = a_1 \ldots a_n$ be a loop in $\G$.  Then
$A(a_i,a_{i+1}) = 1$ for all $1 \leq i \leq n-1$ and
$A(a_n,a_1)=1$.  Hence there exists a loop $\alpha = \alpha_1
\ldots \alpha_n$ in $\Gr(A)$ with $s(\alpha_i) = a_i$ for all
$i$.  Let $f$ be an exit for $\alpha$ in $\Gr (A)$, and
without loss of generality assume $s(f)=s(\alpha_1)$ and $f
\neq \alpha_1$.  Let $b := r(f)$.  Since $A$ has entries in
$\{0,1\}$ we know that $b \neq r(\alpha_1)=a_2$.  Hence $b$
is an exit for $a = a_1 \ldots a_n$.
\end{proof}

\begin{corollary}
Let $\G_1$ and $\G_2$ be ultragraphs with no sinks and
with the same edge matrix $A$.  Then $\G_1$ satisfies
Condition~(L) if and only if $\G_2$ satisfies Condition~(L).
\end{corollary}
\begin{proof} $\G_1$ satisfies Condition~(L)
$\Longleftrightarrow$ $\Gr(A)$ satisfies Condition~(L)
$\Longleftrightarrow$ $\G_2$ satisfies Condition~(L).
\end{proof}

\section{Singular vertices}
\label{sing-sec}
In this section we deal with singular vertices in a
manner similar to what was done in \cite{BPRS} for sinks in
graphs and in \cite{DT} for infinite emitters in graphs.  

\begin{lemma}  Let $\G$ be an ultragraph, let $A$ be a
$C^*$-algebra generated by a Cuntz-Krieger $\G$-family $\{
s_e,p_A \}$, and let $\{ q_n \}$ be a sequence of
projections in $A$.  If $q_n s_\alpha p_A s_\beta^*$
converges for all $\alpha, \beta \in \G^*, A \in \G^0$, then
$\{ q_n \}$ converges strictly to a projection $q \in \M
(A)$.
\label{convergence}
\end{lemma}

\begin{proof}  Since we can approximate any $a \in
A:=C^*(s_e,p_A)$ by a linear combination of $s_\alpha p_A
s_\beta^*$, an $\epsilon / 3$ argument shows that $\{q_n a
\}$ is Cauchy for every $a \in A$.  We define $q : A
\rightarrow A$ by $q(a) := \lim_{n \rightarrow \infty}
q_na$.  Since $$b^* q(a) = \lim_{n \rightarrow \infty}
b^*q_na = \lim_{n \rightarrow \infty} (q_nb)^*a = q(b)^*a,$$
the map $q$ is an adjointable operator on the Hilbert
$C^*$-module $A_A$, and hence defines (left multiplication
by) a multiplier $q$ of $A$ \cite[Theorem 2.47]{RW}. 
Taking adjoints shows that $aq_n \rightarrow aq$ for all $a
\in A$ so $q_n \rightarrow q$ strictly.  It is easy to check
that $q^2 = q = q^*$.
\end{proof}

By adding a tail at a sink $w$ we mean adding a graph of the
form
$$\xymatrix{
w \ar[r]^{e_1} & v_1 \ar[r]^{e_2} & v_2 \ar[r]^{e_3} & v_3
\ar[r]^{e_4} & \cdots}$$
to $\G$ to form a new ultragraph $\F$; thus $F^0 := G^0
\cup \{ v_i : 1 \leq i < \infty \}$, $\F^1 := \G^1 \cup  \{
e_i : 1 \leq i < \infty \}$, and $r$ and $s$ are extended to
$\F^1$ by $r(e_i) = \{ v_i \}$ and $s(e_i) = v_{i-1}$ and
$s(e_1)=w$.  Just as with graphs, when we add tails to sinks
in $\G$ any Cuntz-Krieger $\F$-family will restrict to
a Cuntz-Krieger $\G$-family.  This is because $\F^0$ is
generated by $G^0 \cup \{ v_i : 1 \leq i < \infty \} \cup
\{r(e) : e \in \G^1 \}$ and thus by Lemma~\ref{description}
we see that $\F^0 = \{ A \cup F : A \in \G^0 \text{ and } F
\text{ is a finite subset of } \{v_i \}_{i=1}^\infty \}$.

\begin{proposition}
Let $\G$ be a directed graph and let $\F$ be the ultragraph formed by adding a tail to each sink of $\G$.
\begin{enumerate}
\item For each Cuntz-Krieger $\G$-family $\{S_e, P_A \}$ on
a Hilbert space $\Hi_\G$, there is a Hilbert space $\Hi_\F
:= \Hi_\G \oplus \Hi_T$ and a Cuntz-Krieger $\F$-family
$\{ T_e, Q_A \}$ such that $T_e = S_e$ for $e \in \G^1$,
$Q_A = P_A$ for $A \in \G^0$, and $\sum_{v \notin G^0} Q_{v}$
is the projection on $\Hi_T$. \label{tail1}
\item  If $\{ T_e, Q_A \}$ is a Cuntz-Krieger $\F$-family,
then $\{T_e,Q_A : e \in \G^1, A \in \G^0 \}$ is a
Cuntz-Krieger $\G$-family.  If $w$ is a sink in $E$ such
that $Q_w \neq 0$, then $Q_v \neq 0$ for every vertex on the
tail attached to $w$. \label{tail2}
\item If $\{ t_e, q_A \}$ are the canonical generators of
$C^*(\F)$, then the homomorphism $\pi_{t,q}$ corresponding
to the Cuntz-Krieger $\G$-family $\{t_e,q_A : e \in \G^1, A
\in \G^0 \}$ is an isomorphism of $C^*(E)$ onto a full
corner in $C^*(\F)$. \label{tail3}
\end{enumerate}
\label{beefup}
\end{proposition}

\begin{proof}  For the sake of simplicity we consider the
case in which a single tail has been added to a sink $w$.  As
mentioned earlier, any element $B \in \F^0$ may be
(uniquely) written as $B =A \cup F$ for some $A \in \G^0$
and some finite set $F \subseteq \{ v_i \}_{i=1}^\infty$. 
To extend $\{ S_e, P_A \}$, we let $\Hi_T$ be the direct
sum of infinitely many copies of $P_w \Hi_\G$, define
$P_{v_i}$ to be the projection on the $i^{\text{th}}$
summand, and let $S_{e_i}$ be the partial isometry whose
initial space is the $i^{\text{th}}$ summand and whose final
space is the $(i-1)^{\text{st}}$, with $S_{e_1}$ taking the
first summand of $\Hi_T$ onto $P_w \Hi_\G \subseteq
\Hi_{\G}$.  Now for any $B \in \G^0$ we write $B$ (uniquely)
as $A \cup F$ and define $$P_B := P_A + \sum_{v \in F}
P_v.$$  One can check that $\{ S_e, P_B \}$ is a
Cuntz-Krieger $\F$-family, and hence (\ref{tail1}) holds.

For the same reasons, throwing away the extra elements of a
Cuntz-Krieger $\F$-family gives a Cuntz-Krieger
$\G$-family.  The last statement in (\ref{tail1}) holds
because $$S_{e_1}S_{e_1}^* = P_w \neq 0 \Rightarrow
S_{e_2}S_{e_2}^* = P_{v_1} = S_{e_1}^*S_{e_1} \neq 0 
\Rightarrow S_{e_3}S_{e_3}^* = P_{v_2} = S_{e_2}^*S_{e_2}
\neq 0 \ldots$$
For the first part of (\ref{tail3}), just use part
(\ref{tail1}) to see that every representation of $C^*(\G)$
factors through a representation of $C^*(\F)$. 

We still have to show that the image of $C^*(\G)$ is a full
corner.  List the elements of $G^0 = \{ w_1, w_2, w_3,
\ldots \}$ and the elements $\G^1 = \{ e_1, e_2, e_3, \ldots
\}$.  Define $A_n := \{ v_i : 1 \leq i \leq n \} \cup
\bigcup_{i=1}^n r(e_i)$.  Then given any $e \in \G^1$ and $A
\in \G^0$ we see that for large enough
$n$ we have $p_{A_n} s_e = s_e$ and $p_{A_n} p_A = p_A$. 
Hence Lemma~\ref{convergence} applies and the sequence $\{
p_{A_n}
\}$ converges strictly to a projection $p$ in $\M (C^*(\F))$
satisfying
$$ps_e := \begin{cases} s_e & \text{if $s(e) \in
G^0$} \\ 0 & \text{otherwise}\end{cases} \quad \quad
\text{and} \quad \quad pp_A = p_{A \cap G^0}.$$

\noindent Thus the corner $pC^*(\F)p$ is precisely $C^*(\G)$.

To see that this corner is full suppose $J$ is an ideal
containing $pC^*(\F)p$.  Then $J$ contains $\{q_{r(e)} : e
\in \G^1 \}$ and $\{ q_v : v \in G^0 \}$.  Furthermore, if
$v$ is a vertex on the tail attached to $w$, then there is a
unique path $\alpha$ with $s(\alpha) = w$ and $r(\alpha)=v$,
and $$q_w \in J \Longrightarrow t_\alpha = q_w t_\alpha \in J
\Longrightarrow q_v = t_\alpha^* t_\alpha \in J.$$  Thus $J$
contains $\{ q_{r(e)} : e \in \F^1 \} \cup \{q_v : v \in F^0
\}$ and hence is all of $C^*(\F)$.
\end{proof}

Now suppose  that $\G$ is an ultragraph with an infinite
emitter $v_0$.  We \emph{add a tail} at $v_0$ by
performing the following procedure.  List the edges $g_1,
g_2, \ldots$ of $s^{-1}(v_0)$.  We begin by adding
vertices and edges as we did with sinks:
$$
\xymatrix{
v_0 \ar[r]^{e_1} & v_1 \ar[r]^{e_2} & v_2 \ar[r]^{e_3} & v_3
\ar[r]^{e_4} & \cdots } $$
Then we remove the edges in $s^{-1}(v_0)$ from $\F$, and for
each $j$ we draw an edge $f_j$ with source $v_{j-1}$ and
range $r(g_j)$.   

To be precise, if $\G$ is an ultragraph with an infinite
emitter $v_0$, we define $F^0 := G^0 \cup \{v_1, v_2,
\dots\}$ and 
$$\F^1 := \{e \in G^1 : s(e) \neq v_0 \} \cup
\{e_i\}_1^\infty \cup
 \{ f_j : 1 \leq j < \infty \}.$$
We extend $r$ and $s$ to $\F^1$ as indicated above.  In
particular, $s(e_i) = v_{i-1}$, $r(e_i) = \{ v_i \}$,
$s(f_j) = v_{j-1}$, and $r(f_j) = r(g_j)$.

For any $j$ we shall often have need to refer to the path
$\alpha^j := e_1e_2\ldots e_{j-1} f_j$ in $\F$.  Also note
that if $\G$ is an ultragraph and $\F$ is the graph formed
by adding a tail at an infinite emitter, then $\F^0$ is
generated by $\{ \{v \} : v \in F^0 \} \cup \{r(e) : e \in
\F^1 \} = \{ \{v \} : v \in G^0 \} \cup \{ r(e) : e \in \G^1
\} \cup \{ \{v_i \} : 1 \leq i < \infty \}$.  Thus by
Lemma~\ref{description} we see that $\F^0 = \{ A \cup F : A
\in \G^0 \text{ and } F \text{ is a finite subset of } \{v_i
\}_{i=1}^\infty \}$. 

\begin{definition}
If $\G$ is an ultragraph, a \emph{desingularization of
$\G$} is an ultragraph $\F$ obtained by adding a tail at
every singular vertex of $\G$. 
\end{definition}

\begin{lemma}
\label{lem-CKEinF}
Suppose $\G$ is a graph and let $\F$ be a
desingularization of $\G$.  
If $\{T_e, Q_A\}$ is a Cuntz-Krieger $\F$-family, then there
exists a  Cuntz-Krieger $\G$-family in $C^*(\{T_e, Q_A\})$.  
\end{lemma}

\begin{proof} For every $A \in \G^0$, define
$P_A := Q_A$.  For every edge $e \in \G^1$ with $s(e)$ not a
singular vertex, define $S_e := T_e$.  If $e \in \G^1$ with
$s(e) = v_0$, a singular vertex, then $e = g_j$ for some
$j$ and we define $S_e := T_{\alpha^j}$.   The fact that
$\{S_e, P_A : e \in \G^1, A \in \G^0\}$ is a Cuntz-Krieger
$\G$-family follows immediately from the fact that $\{T_e,
Q_A : e \in \F^1, A \in \F^0 \}$ is a Cuntz-Krieger
$\F$-family.
\end{proof}

\begin{lemma}
\label{lem-extendEfamtoF}
Let $\G$ be an ultragraph and let $\F$ be 
a desingularization of $\G$. For every
Cuntz-Krieger $\G$-family $\{S_e, P_A : e \in \G^1, A \in
\G^0\}$ on a Hilbert space $\Hi_\G$, there exists a
Hilbert space $\Hi_\F = \Hi_\G \oplus
\Hi_T$ and a   Cuntz-Krieger $\F$-family $\{T_e, Q_A
: e \in \F^1,  A \in \F^0\}$ on $\Hi_\F$ satisfying: 
\begin{itemize}
\item $P_A = Q_A$ for every $A \in \G^0$;
\item $S_e = T_e$ for every $e \in \G^1$ such that $s(e)$ is
not a singular vertex;
\item $S_e = T_{\alpha^j}$ for every $e = g_j \in \G^1$
such that $s(g_j)$ is a singular vertex;
\item $\sum_{v \notin G^0} Q_v$ is the projection
onto $\Hi_T$.
\end{itemize}
\end{lemma}

\begin{proof}
We prove the case where $\G$ has just one singular vertex
$v_0$. If $v_0$ is a sink, then the result follows from
Proposition~\ref{beefup}.  Thus we need only consider
when $v_0$ is an infinite emitter.  Given a Cuntz-Krieger
$\G$-family $\{S_e, P_A\}$, and a nonnegative integer $n$ we
define $R_0 = 0$ and $R_n := \sum_{j=1}^n S_{g_j}S_{g_j}^*$. 
Note that the $R_n$'s are projections because the
$S_{g_j}$'s have orthogonal ranges.  Furthermore, $R_n \leq
R_{n+1} < P_{v_0}$ for every $n$.

Now, for every integer $n \geq 1$, define $\Hi_n
:= (P_{v_0} - R_n)\mathcal{H}_\G$ and set  
$$ \Hi_\F := \Hi_\G \oplus \bigoplus_{n=1}^{\infty}
\Hi_n. $$

\noindent As mentioned previously any $B \in \F^0$ may
be written (uniquely) as $A \cup F$ for some $A \in \G^0$
and some finite set $F$ of vertices on the added tail.

For every $A \in \G^0$,
define $Q_A$ to equal $P_A$ on the $\Hi_\G$ component
of $\Hi_\F$ and zero elsewhere.  That is, $Q_A (\xi_\G,
\xi_1, \xi_2, \ldots) = (P_A \xi_\G, 0, 0, \ldots)$. 
Similarly, for every $e \in \G^1$ with
$s(e) \neq v_0$, define $T_e = S_e$ on the $\Hi_\G$
component; $T_e (\xi_\G, \xi_i, \xi_2, \ldots) = (S_e
\xi_\G, 0, 0, \ldots )$.  For each vertex
$v_n$ on the added tail, define $Q_{v_n}$ to be the
projection onto $\Hi_n$; $Q_{v_n}(\xi_\G,
\xi_1, \ldots ,\xi_n, \xi_{n+1}, \ldots) =
(0,0,\ldots,\xi_n,0,\ldots)$.   Note that, because the
$R_n$'s are non-decreasing,
$\Hi_{n} \subseteq \Hi_{n-1}$ for each $n$.  Thus,
for each edge $e_n$ of the tail, we may define $T_{e_n}$ to
be the inclusion of $\Hi_n$ into
$\Hi_{n-1}$ (where $\Hi_0$ is taken to mean
$P_{v_0}(\Hi_\G$)).  More precisely, 
$$ T_{e_n}(\xi_\G, \xi_1, \xi_2, \ldots) =
(0,0,\ldots,0,\xi_n,0,\ldots), $$
where the $\xi_n$ is in the $\Hi_{n-1}$ component.  

Finally, for each edge $g_j$ and for each $\xi \in
\Hi_\G$, we have that $S_{g_j}\xi \in
\Hi_{j-1}$.  Thus, we can define $T_{f_j}$ as
$$ T_{f_j}(\xi_\G, \xi_1, \xi_2, \ldots) =
(0,\ldots,0,S_{g_j}\xi_\G,0,\dots), $$
where the nonzero term appears in the $\Hi_{j-1}$ component.

Now for any $B \in \F^0$ we can (uniquely) write $B := A \cup
F$ for some $A \in \G^0$ and some finite subset $F$ of
vertices on the added tail.  Thus we may define $Q_B := Q_A
+ \sum_{v \in F} Q_v$.  It then follows from
calculations similar to those in \cite[Lemma 2.10]{DT} that
$\{ T_e, Q_A \}$ is a Cuntz-Krieger $\F$-family satisfying
the bulleted points.
\end{proof}

\begin{proposition}
\label{desingularization}
Let $\G$ be an ultragraph and let $\F$ be a
desingularization of $\G$.  Then $C^*(\G)$ is isomorphic to a
full corner of $C^*(\F)$.
\end{proposition}

\begin{proof}
Again, for simplicity we assume $\G$ has only one singular
vertex $v_0$.  If
$v_0$ is a sink, then the claim follows from
Proposition~\ref{beefup}.  Therefore, let us assume that
this singular vertex is an infinite emitter.  Let
$\{t_e, q_A : e \in \F^1, A \in \F^0\}$ denote the
canonical set of generators for $C^*(\F)$ and let $\{s_e, p_A
: e \in \G^1, A \in \G^0\}$ denote the Cuntz-Krieger
$\G$-family in $C^*(\F)$ constructed in
Lemma~\ref{lem-CKEinF}.  Define $B:=C^*(\{s_e, p_A\})$. 
Also, list the elements of $G^0 = \{ w_1, w_2, \ldots \}$
and the elements of $\G^1 = \{ h_1, h_2, \ldots \}$.   For
each nonnegative integer $n$ let $A_n := \{w_i : 1 \leq i
\leq n \} \cup \bigcup_{i=1}^n r(h_i)$.  It follows
from Lemma~\ref{convergence} that the sequence $\{ p_{A_n}
\}$ converges to a projection $p \in \M (C^*(\F))$ satisfying
$$pt_e := \begin{cases} t_e & \text{if $s(e) \in G^0$} \\ 0
& \text{otherwise}\end{cases} \quad \quad
\text{and} \quad \quad pq_A = q_{A \cap G^0}.$$  From these
relations one can see that $B \cong pC^*(\F)p$.

We shall now show that $B \cong C^*(\G)$.  Since $B$ is
generated by a Cuntz-Krieger $\G$-family, it suffices to show
that $B$ satisfies the universal property for
$C^*(\G)$.  Let $\{S_e, P_A : e \in \G^1,
A \in \G^0 \}$ be a Cuntz-Krieger $\G$-family on a Hilbert
space $\Hi_\G$.  Then by Lemma~\ref{lem-extendEfamtoF}
we can construct a Hilbert space $\Hi_\F$ and a
Cuntz-Krieger $\F$-family $\{T_e, Q_A : e \in \F^1, A \in
\F^0 \}$ on $\Hi_\F$ such that $Q_A = P_A$ for every $A \in
\G^0$, $T_e = S_e$ for every $e \in \F^1$ with $s(e) \neq
v_0$, and $S_{g_j} = T_{\alpha^j}$ for every edge
$g_j \in \G^1$ whose source is $v_0$.  By the universal
property of $C^*(\F)$ we have a homomorphism
$\pi$ from $C^*(\F)$ onto $C^*(\{T_e, Q_A : e \in \F^1, A
\in \F^0 \})$ that takes
$t_e$ to $T_e$ and $q_A$ to $Q_A$.   Now $p_A = q_A$ for any
$A \in \G^0$, so $\pi(p_A) = Q_A = P_A$.  Let $e
\in \G^1$ and $s(e) \neq v_0$.  Then $s_e = t_e$ and
$\pi(s_e) = T_e = S_e$.  Finally, if $s(e) = v_0$, then $e =
g_j$ for some $j$, $s_e = t_{\alpha^j}$, and
$\pi(s_{g_j}) = T_{\alpha^j} = S_{g_j}$.
Thus $\pi|_B$ is a representation of $B$ on $\Hi_\G$ that
takes generators of $B$ to the corresponding elements of
the given Cuntz-Krieger $\G$-family.  Therefore $B$ satisfies
the universal property of $C^*(\G)$ and $C^*(\G) \cong B$.

Finally, we note that the corner $C^*(\G) \cong B
\cong pC^*(F)p$ is full by an argument similar to the one
given in Proposition~\ref{beefup}.
\end{proof}

These results give us a way to extend the uniqueness
theorems of \S\ref{uniqueness-sec} to
$C^*$-algebras of ultragraphs that contain singular
vertices.  Also note that an ultragraph satisfies
Condition~(L) if and only if its desingularization satisfies
Condition~(L).

\begin{theorem}[Uniqueness]
Suppose that $\G$ is an ultragraph and that every loop in
$\G$ has an exit.  If $\{S_e,P_A \}$ and $\{T_e, Q_A \}$ are
two Cuntz-Krieger $\G$-families in which all the projections
$P_A$ and $Q_A$ are nonzero, then there is an isomorphism
$\phi$ of $C^*(S_e,P_A)$ onto $C^*(T_e,Q_v)$ such that
$\phi(S_e)=T_e$ and $\phi(P_A) = Q_A$.
\label{uniqueness}
\end{theorem}

\begin{proof}  Let $\F$ be a desingularization of $\G$. 
Then we may use Lemma~\ref{lem-extendEfamtoF} to
extend the $\G$-families to $\F$-families.  Since every loop
in $\G$ has an exit, it follows that every loop in $\F$ has
an exit, and thus we may apply
Proposition~\ref{uniqueness-nosing} to get an isomorphism
$\phi$ that restricts to our desired isomorphism from
$C^*(S_e,P_A)$ onto $C^*(T_e,Q_v)$.
\end{proof}

\begin{theorem}[Gauge-Invariant Uniqueness]
Let $\G$ be an ultragraph, $\{s_e,p_A\}$ the canonical
generators in $C^*(\G)$, and $\gamma$ the gauge action on
$C^*(\G)$.  Also let $B$ be a $C^*$-algebra, and $\phi :
C^*(\G) \rightarrow B$ be a homomorphism for which
$\phi(p_A) \neq 0$ for all nonempty $A$.  If there
exists a strongly continuous action $\beta$ of $\T$ on $B$
such that $\beta_z \circ \phi = \phi \circ \gamma_z$ for all
$z \in \T$, then $\phi$ is faithful.
\label{g-i-uniqueness}
\end{theorem}

\begin{proof}
Since $\beta : \T \rightarrow \aut B$ is an action of $\T$
on $B$, there exists a Hilbert space $\Hi_\G$, a faithful
representation $\pi : B \rightarrow \B (\Hi_\G)$, and a
unitary representation $U : \T \rightarrow U (\Hi_\G)$ such
that $$\pi(\beta_z(x)) = U_z \pi(x) U_z^* \quad \quad \text{
for all $x \in B$ and $z \in \T$.}$$  Let $S_e := \pi \circ
\phi (s_e)$ and $P_A := \pi \circ \phi (p_A)$.  Also let $\F$
be a desingularization of $\G$.  For simplicity, we shall
assume that $\G$ has only one singular vertex $v_0$.

If $v_0$ is a sink, then it follows from
Proposition~\ref{beefup} that there exists a Hilbert space
$\Hi_\F := \Hi_\G \oplus \Hi_T = \Hi_\G \oplus
\bigoplus_{i=1}^\infty P_{v_0} \Hi_\G$, and a Cuntz-Krieger
$\F$-family $\{T_e, Q_A \}$ in
$\B (\Hi_\F)$ that restricts to $\{ S_e, P_A \}$.  We shall
define a unitary representation $V : \T \rightarrow
U(\Hi_\F)$ as follows:  If $h \in \Hi_\G$, then we define
$V_z h := U_z h \in \Hi_\G$.  If $h \in Q_{v_i} \Hi_\G =
P_{v_0} \Hi_\G$ is in the $i^{\text{th}}$ component of
$\Hi_T$, then we define $V_z h := z^{-i} U_z h \in Q_{v_i}
\Hi_\G$ also in the $i^{\text{th}}$ component of $\Hi_T$.  We
define $V_z$ on all of $\Hi_\F$ by extending it linearly.

Now let $\pi_{t,q} : C^*(\F) \rightarrow \B (\Hi_\F)$,
be the representation for which $\pi_{t,q} (t_e) = T_e$ and
$\pi_{t,q}(q_A) = Q_A$, and let $\gamma^\F$ be the gauge
action on $C^*(\F)$.  Then one can check that $$\pi_{t,q}
(\gamma^\F_z(x)) = V_z \pi_{t,q}(x) V_z^* \quad \quad
\text{for all $x \in C^*(\F)$ and all $z \in \T$.}$$  (To see
this simply check the relation on the generators $\{t_e \}$
and use the fact that $\pi \circ \phi (s_e) = \pi_{t,q}
(t_e)$ for all $e \in \G^1$.)  Now if we define
$\tilde{\beta} : \T \rightarrow \aut C^*(T_e,Q_A)$ by
$\tilde{\beta}_z(X) := V_z X V_z^*$, then we see that
$\tilde{\beta}_z \circ \pi_{t,q} =
\pi_{t,q} \circ \gamma^F_z$ for all $z \in \T$.  Since
$\pi_{t,q}(q_A) = Q_A \neq 0$ and since $\F$ has no sinks or
infinite emitters, it follows from
Proposition~\ref{gauge-invariant-nosing} that $\pi_{t,q}$ is
faithful.  Now, if $\iota : C^*(\G) \rightarrow C^*(\F)$
denotes the canonical inclusion of $C^*(\G)$ into $C^*(\F)$,
then we see that $\pi_{t,q} \circ \iota = \pi \circ \phi$
(since each map agrees on the generators $\{s_e,p_A \}$). 
Because $\iota$ and $\pi_{t,q}$ are both injective, it
follows that $\phi$ is injective.

If $v_0$ is an infinite emitter, then an argument almost
identical to the one above works.  We obtain a faithful
representation $\pi : C^*(\G) \rightarrow \Hi_\G$ and a
unitary representation $U : \T \rightarrow U (\Hi_\G)$ as
before, and we then extend this Hilbert space to $\Hi_\F :=
\Hi_\G \oplus \bigoplus_{n=1}^\infty \Hi_n$, where $\Hi_n :=
(P_{v_0} - \sum_{j=1}^{n-1} S_{g_j}S_{g_j}^*)(\Hi_\G)$ as
in Lemma~\ref{lem-extendEfamtoF}.  Similarly, we define $V :
\T \rightarrow \B (\Hi_\F)$ as follows:  If $h \in \Hi_\G$,
then $V_z h := U_z h \in \Hi_\G$.  If $h \in \Hi_n \subseteq
\Hi_\G$, then $V_z h := z^{-n} U_z h \in \Hi_n$.  The rest
of the argument follows much like the one above.

If $\G$ has more than one sink or more than one infinite
emitter, then we simply account for multiple tails.  The
previous argument will still work, we need only keep track
of the multiple pieces added on when extending $\Hi_\G$ to
$\Hi_\F$.
\end{proof}

\section{Concluding Remarks}

We have seen in this paper that ultragraph algebras
contain both the Exel-Laca algebras and the graph
algebras.  Furthermore, it is shown in a forthcoming article
\cite{Tom4} that there exist ultragraph algebras that
are neither Exel-Laca algebras nor graph algebras. 
Throughout this paper we have seen that many of the
techniques of graph algebras can be applied to ultragraph
algebras and that analogues of the results for graph algebras
and Exel-Laca algebras hold for ultragraph algebras. 

These observations are important for many reasons.  First of
all, ultragraphs give a context in which many results
concerning graph algebras and Exel-Laca
algebras may be proven simultaneously.  In the past, many
similar results (e.g. the Cuntz-Krieger Uniqueness Theorem,
the Gauge-Invariant Uniqueness Theorem) were proven
separately for graph algebras and for Exel-Laca algebras. 
Since ultragraph algebras contain both the graph algebras
and the Exel-Laca algebras, it suffices to prove these
results once for ultragraph algebras.  Hence these
classes are in some sense unified under the umbrella of
ultragraph algebras.  Second, since many of the
graph techniques may be used for ultragraphs, we
see that we may study Exel-Laca algebras in this context and
the (often complicated) matrix techniques may be avoided in
favor of graph techniques.  Finally, ultragraph
algebras are a larger class of
$C^*$-algebras than the Exel-Laca algebras and the graph
algebras.  Thus with only slightly more work, we are able to
extend these results to a larger class of $C^*$-algebras.

\end{document}